\documentclass{amsart}
\usepackage{amssymb}
\usepackage{amsfonts}
\usepackage{amsmath}
\usepackage[legalpaper,bookmarks=true,colorlinks=true,linkcolor=blue,citecolor=blue]{hyperref}
\usepackage{graphicx}%
\usepackage{fancyhdr}
\usepackage{color}
\usepackage[mathlines]{lineno}
\usepackage{lscape}
\usepackage{epsfig}
\usepackage{natbib}
\usepackage{geometry}
\usepackage{tgbonum}
\fontfamily{qcr}\selectfont

\theoremstyle{plain}

\newtheorem{fact}{Fact}

\numberwithin{equation}{section}

\newcommand{\Bin}{\bigskip \noindent}

\newcommand{\Ni}{\noindent}

\begin{document}
\Large
\title[]{A direct construction of the standard Brownian motion}

\begin{abstract} In this note, we combine the two approaches of Billingsley (1998) and Cs\H{o}rg\H{o} and R\'ev\'esz (1980), to provide a detailed sequential and descriptive for creating s standard Brownian motion, from a  Brownian motion whose time space is the class of non-negative dyadic numbers. By adding the proof of Etemadi's inequality to text, it becomes self-readable and serves as an independent source for researches and professors.\\  

\noindent $^{\dag}$ Gane Samb Lo.\\
LERSTAD, Gaston Berger University, Saint-Louis, S\'en\'egal (main affiliation).\newline
LSTA, Pierre and Marie Curie University, Paris VI, France.\newline
AUST - African University of Sciences and Technology, Abuja, Nigeria\\
gane-samb.lo@edu.ugb.sn, gslo@aust.edu.ng, ganesamblo@ganesamblo.net\\
Permanent address : 1178 Evanston Dr NW T3P 0J9,Calgary, Alberta, Canada.\\

\noindent Aladji Babacar Niang\\
LERSTAD, Gaston Berger University, Saint-Louis, S\'en\'egal.\\
Email: aladjibacar93@gmail.com\\

\noindent Dr Harouna SANGARE\\
Main Affiliation: DER MI, FST, Universit\'e des Sciences, des Techniques et des Technologies de Bamako (USTT-B), Mali.\\
Affiliation : LERSTAD, Universit\'e Gaston Berger (UGB), Saint-Louis, S\'en\'egal.\\
Email : harounasangare@fst-usttb-edu.ml, harouna.sangare@mesrs.ml\\
sangare.harouna@ugb.edu.sn, harounasangareusttb@gmail.com\\

\end{abstract}
\maketitle

\Bin \textbf{NB} This work has been completed while the authors were in confinement in Saint-Louis (SENEGAL), due to the COVID-19 pandemy.\\

\section{introduction} \label{cmb_01_sec_01}

\Ni In many textbooks, if not the bulk of them, the Kolmogorov Existence Theorem (\textit{KET}) (see \cite{ips-mfpt-ang}, Chapter 10) is used to directly get a stochastic process $(B_s)_{s \in \mathbb{R}_{+}}$ having the probability law of a Brownian motion, that is a centering and Gaussian stochastic process such that for all $(s,t) \in [0,+\infty[^2$, $\mathbb{E}(B_s B_t)=min(s,t)$. Afterward, some elements of random analysis are used to transform that process to an \textit{a.s.} continuous process called Standard Brownian motion. The reader is referred to \cite{billinsgleypm} for example for that approach.\\

\Ni Other authors, \cite{csorgoR} for example, use a more constructive method. In \cite{csorgoR}, in a first step, the Kolmogorov Existence Theorem (\textit{KET}) is used to create a probability space holding a family $(X_s)_{s \in D}$ of independent $\mathcal{N}(0,1)$-random variables, indexed by the space time $D$ of non-negative dyadic  numbers. In a second step, a stochastic process $(B_s)_{s \in D}$ is constructed from  $(X_s)_{s \in D}$, which should follow the law of a Brownian motion of $D^{+}$. A third step should be an extension of  $(B_s)_{s \in D}$ to a standard Brownian motion $(B_s)_{s \in \mathbb{R}_+}$.\\

\Ni The second construction has particular advantages when it comes to study the extrema of Brownian processes and in approximation theory of empirical functions by Brownian bridges. For that reason, we wish to devote a fully documented note of the second construction. Beside, that not would serve as an innovative pedagogical document.\\

\Ni To highlight some contribution of that note, we should remark that :\\

\Ni (a) \cite{csorgoR} justified the existence of the Standard Brownian continuous extension and indicate that checking it has a Brownian law is easy to get.\\

\Ni (b) The argument used in the extension are not proved in that document and references on their are not given.\\

\Ni (c) Fortunately, the arguments in \cite{billinsgleypm} are enough to realize the needed extension.\\

\Ni In this paper, we combine the two texts and provide a fully documented note can be used for graduate students, professors and researchers as an independent document.\\

\Ni Let us begin by the first step.

\section[A sequential construction of a Brownian motion]{A sequential construction of a Brownian motion of the set of non-negative dyadic  numbers} \label{cmb_02_sec_01} \label{cmb_02_sec_01}

\Ni Let us denote

$$
D=\left\{\frac{k}{2^n}, \\ k \in \mathbb{N}, \ \ n\geq 0 \right\}. 
$$

\Bin We begin by an easy application of the \textit{KET}: there exists a probability space $\left(\Omega,\mathcal{A},\mathbb{P}\right)$ bearing a family $\mathcal{F}=(X_r)_{r\in D}$ of independent $\mathcal{N}(0,1)$ random variables. The justification is given in \cite{ips-mfpt-ang}, Chapter 9, Section 6, Example 1.\\

\Bin We recall that $\mathbb{N}$ is the class of all non-negative integers, including $0$ and $\mathbb{N}^{*}=:\mathbb{N}\setminus \{0\}$ is the class of positive integers. From that family $\mathcal{F}$, we create a stochastic process by induction as follows: \\

\Ni (i) $B(0) = 0$; \\

\Ni (ii) $\forall k\in \mathbb{N}^{*}, B(k) = X_1 + \cdots + X_k$; \\ 

\Ni (iii) $\forall \ell\in \mathbb{N}$, 

\begin{equation}
B\left(\frac{2\ell+1}{2}\right) = \frac{B(\ell)+B(\ell+1)}{2} + \frac{X_{(2\ell+1)/2}}{\sqrt{2^2}}. \label{ST2}
\end{equation}

\Bin Let us explain that step in details. Here, $n=1$ and we have for all $k\geq 1$,

\begin{equation*}
\frac{k}{2} = \frac{2\ell+s}{2}, \ s\in \left\{0,1\right\}.
\end{equation*}

\Bin If $s=0$, $B(\frac{k}{2})=B(\ell)$ is already defined in Step (i) or Step (ii). If $s=1$, Formula \label{ST2} becomes 

\begin{eqnarray} \label{defStep2}
B\left(\frac{2\ell+1}{2^n}\right) &=& \frac{B(\frac{2\ell}{2^{n}})+B(\frac{2\ell+2}{2^{n}})}{2} + \frac{X_{(2\ell+1)2^{-n}}}{\sqrt{2^{n+1}}}\notag\\
&=& \frac{B(\frac{\ell}{2^{n-1}})+B(\frac{\ell+1}{2^{n-1}})}{2} + \frac{X_{(2\ell+1)2^{-n}}}{\sqrt{2^{n+1}}}.
\end{eqnarray}

\Bin That principle we applied from $n=0$ to $n=1$ may be repeated from $n=1$ to $n=2$, etc. in the following way: \\

\Ni Given 

$$
\left\{B\left(\frac{k}{2^j}\right), \ k\in \mathbb{N}, \ j\geq 0\right\}
$$ 

\Bin is defined up to $j=n$, we pass to step $n+1$ by re-conducting Formula \eqref{defStep2} as follows

\begin{equation} \label{defStepnpun}
B\left(\frac{2\ell+1}{2^{n+1}}\right) = \frac{B(\frac{\ell}{2^{n}})+B(\frac{\ell+1}{2^{n}})}{2} + \frac{X_{(2\ell+1)2^{-(n+1)}}}{\sqrt{2^{n+2}}}.
\end{equation}

\Bin The stochastic process 

$$
B=\left\{B(r), \ r\in D\right\}
$$

\Bin is entirely defined. \\

\Ni Let us show that $B$ satisfies the conditions of Brownian movement on the time space $D$, \textit{i.e.} \\

\Ni (IC) $B(0)=0$; \\

\Ni (GM) $\forall r\in D$,\ $B(r)\sim \mathcal{N}(0,r)$; \\

\Ni (II) $\forall k\geq 2$,\ $\forall (r_1, \cdots, r_k)\in D^k$ such that $0=r_0<r_1< \cdots <r_k$, the vector  $\left(B(r_1)-B(r_0),B(r_2)-B(r_1), \cdots , B(r_k)-B(r_{k-1})\right)^t$ has independent margins; \\

\Ni (IC) $\forall 0\leq r_1<r_2$,\ $B(r_2)-B(r_1)=_d B(r_2-r_1)$.\\

\Ni \textbf{Proof}. First we make general remarks and next proceed by induction. Here are the general remarks. Let us denote, for $n\geq 0$,

$$
D_{n}=\left\{\frac{k}{2^{n}}, \ k \in \mathbb{N}\right\} \ and \ D_{n}^{\star}=\left\{\frac{2k+1}{2^{j}}, \ k \in \mathbb{N}, \ 0\leq j \leq n\right\}.
$$ 

\Bin Actually, $D_{n}^{\star}$, for $n\geq 0$, is the class of all dyadic numbers which can be reduced to a form $k/2^{j}$, $k\geq 0$ and $0\leq j\leq n$. In that view, we have for $\ell\geq 0$ and $n\geq 0$,

$$
\frac{2\ell+1}{2^{n+1}} \notin D_{n}^{\star}.
$$

\Bin So we have the following remarks from the process of creating $B$.\\

\Ni (R1) For all $n\geq 0$, for $r \in D_{n}^{\star}$, $B(r)$ is function of quantities  $X_{s}$ with $s \in D_{n}^{\star}$.\\

\Ni (R2) Hence, for all $n\geq 0$ and for $r \in D_{n}^{\star}$, $B(r)$ is a finite linear combination of elements of $(X_{s})_{s\in D_{n}^{\star}}$.\\

\Ni (R3) Hence, all margins $B(r)$, $r \in D$, are finite linear combination for the Gaussian process $(X_r)_{r \in D}$ which has independent 
$\mathcal{N}(0,1)$-components. So $B$ is Gaussian on $D$.\\

\Ni (R4) By defining $B((2\ell+1)/2^{n+1})$ in \eqref{defStepnpun}, the quantities

$$
R(1,\ell,n+1)=:\frac{B(\frac{\ell}{2^{n}})+B(\frac{\ell+1}{2^{n}})}{2}  \ and \ R(2,\ell,n+1)=\frac{X_{(2\ell+1)2^{-(n+1)}}}{\sqrt{2^{n+2}}}
$$

\Bin are independent.\\

\Ni (R5) So Checking condition (II) is based on checking that the co-variance of the element of that vector are zero. So, it is enough to check it for $k=2$.\\

\Ni Now, we may continue the proof by induction. Based on the remarks above, and only the fact that on a finite number of the indices (r) is used in checking any of the properties (GM), (II) and (IC), we can see that is enough to establish, for any $n\geq 1$,  \\

\Ni (GM) $\forall r\in D_{n}$,\ $B(r)\sim \mathcal{N}(0,r)$. \\

\Ni (II)  $\forall (r_1,r_2,r_3,r_4) \in (D_{n}^{\star})^4$ such that $0=r_0<r_1<r_2 \leq r_3< r_4$, $B(r_4)-B(r_3)$ and $B(r_2)-B(r_{1})$  are independent. \\

\Ni (IC) $\forall (r_1,r_2) \in (D_{n}^{\star})^2$ such that $0=r_0<r_1<r_2$,  $B(r_2)-B(r_1)=_d B(r_2-r_1)$.\\

\Bin  Let us proceed by induction.\\

\Ni \textit{Initial condition}. For $n=0$. For $r \in D_{1}^{\star}$, we have either $B(r)=0$ (for $r=0$) or $B(r)=X_1+\cdots+X_k$ (for $r=k$). Hence (GM), (II) and (IS) are obvious.\\

\Ni \textit{Induction hypothesis}. Suppose that (GM), (II) and (IS) hold on $D_{n}^{\star}$. The following fact is standard and left as exercise.

\begin{fact} Let $\emptyset \neq T\subset \mathbb{R}$. Given that $\{X(t), \ t\in T\}$ is a centered and Gaussian stochastic process, the set of the following four conditions :\\

\Ni (IC) $X(0) =0$, a.s.,\\

\Ni (GM) $\forall t \in T$,\ $X(t)\sim \mathcal{N}(0,t)$, \\

\Ni (II)  $\forall (t_1,t_2,t_3,t_4) \in T^4$ such that $t_1<t_2 \leq t_3< t_4$, $X(t_4)-X(t_3)$ and $X(t_2)-X(t_{1})$  are independent, \\

\Ni (IC) $\forall (t_1,t_2) \in T^2$ such that $t_1<t_2$,  $X(t_2)-X(t_1)=_d X(t_2-t_1)$,\\

\noindent is equivalent to

$$
\forall (s,t) \in T^2, \ \mathbb{C}ov(X(t),X(s))=\min(s,t).
$$
\end{fact}
  
\Bin By applying this, the induction assumption implies that

$$
\forall  (r_1,r_2) \in (D_{n}^{\star})^2, \mathbb{C}ov(B(r_1), B(r_2))=\min(r_1,r_2).
$$

\Bin Let us check each of the properties (GM), (II) and (IC) at level $n+1$.\\

\Ni \textit{(a) Checking (GM)}. For $r\in D_{n+1}^{\star}$, we have $r=(2\ell)2^{-(n+1)}$ or $r=(2\ell+1)2^{-(n+1)}$, $\ell\geq 0$. The first case in covered by the induction hypothesis. It remains to do the work for the second case. But by definition, 

$$
B\left(\frac{2\ell+1}{2^{n+1}}\right)=R(1,\ell, n+1) + R(2,\ell, n+1)
$$

\Bin with

$$
R(1,\ell, n+1)=\frac{B(\frac{\ell}{2^{n}})+B(\frac{\ell+1}{2^{n}})}{2}  \ and \ R(2,\ell,n+1)=\frac{X_{(2\ell+1)2^{-(n+1)}}}{\sqrt{2^{n+2}}}
$$

\Bin Since $R(1,\ell, n+1)$ and $R(2,\ell, n+1)$ are centered and independent Gaussian random variables, we have

$$
\mathbb{V}ar\left( B\left(\frac{2\ell+1}{2^{n+1}}\right)\right)= \mathbb{V}ar(R(1,\ell, n+1))+ \mathbb{V}ar(R(2,\ell, n+1)).
$$

\Bin But

$$
\mathbb{V}ar(R(2,\ell, n+1))=\frac{1}{2^{n+2}}
$$

\Bin and, by using

\begin{eqnarray*}
\mathbb{V}ar(R(1,\ell, n+1))&=&\frac{1}{4}\left(\mathbb{V}ar\left( B\left(\frac{\ell}{2^{n}}\right)\right) + \mathbb{V}ar\left( B\left(\frac{\ell+1}{2^{n}}\right)\right) + 2\mathbb{C}ov\left( B\left(\frac{\ell}{2^{n}}\right), B\left(\frac{\ell+1}{2^{n}}\right) \right)\right)\\
&=&\frac{1}{4}\left(\frac{\ell}{2^{n}} + \frac{\ell+1}{2^{n}} + 2\frac{\ell}{2^{n}}\right)\\
&=&\frac{3\ell}{2^{n+2}} + \frac{\ell+1}{2^{n+2}}.
\end{eqnarray*}

\Bin We get

$$
\mathbb{V}ar\left( B\left(\frac{2\ell+1}{2^{n+1}}\right)\right)=\frac{4\ell+2}{2^{n+2}}=\frac{2\ell+1}{2^{n+1}}=r.
$$

\Bin \textit{(b) Checking (II)}. For $(r_1,r_2,r_3,r_4) \in (D_{n+1}^{\star})^4$, such that $0\leq r_1<r_2<r_3<r_4$. Since the process is centered and Gaussian, checking (II) reduces to checking that

$$
\mathbb{C}ov(B(r_2)-B(r_1), B(r_4)-B(r_3))=0.
$$

\Bin But, for each $i \in \{1,2,3,4\}$, we may have either $r_i=(2\ell_i)2^{-(n+1)}$ or $r_i=(2\ell_i+1)2^{n+1}$. So we have to check sixteen cases. But we may exclude the case in which $r_i=(2\ell_i)2^{-(n+1)}$ for all $i \in \{1,2,3,4\}$ since this is covered by the induction hypothesis. Here, we only check the case
$r_i=(2\ell_i+1)2^{-(n+1)}$ for all $i \in \{1,2,3,4\}$ for peace of mind. The treatment of the fourteen other cases are very similar.

\Bin By using the induction hypothesis, the remarks and the definitions

\begin{equation} 
B\left(\frac{2\ell_i+1}{2^{n+1}}\right) = \frac{B(\frac{\ell_i}{2^{n}})+B(\frac{\ell_i+1}{2^{n}})}{2} + \frac{X_{(2\ell_i+1)2^{-(n+1)}}}{\sqrt{2^{n+2}}},
\end{equation}

\Bin with $\ell_1<\ell_2\leq \ell_3<\ell_4$, we can develop the product of

$$
Y1=\frac{B(\frac{\ell_2}{2^{n}})+B(\frac{\ell_2+1}{2^{n}})}{2} + \frac{X_{(2\ell_2+1)2^{-n}}}{\sqrt{2^{n+2}}}-\frac{B(\frac{\ell_1}{2^{n}})-B(\frac{\ell_1+1}{2^{n}})}{2} - \frac{X_{(2\ell_1+1)2^{-n}}}{\sqrt{2^{n+2}}}
$$

\Bin by

$$
Y2=\frac{B(\frac{\ell_4}{2^{n}})+B(\frac{\ell_4+1}{2^{n}})}{2} + \frac{X_{(2\ell_4+1)2^{-(n+1)}}}{\sqrt{2^{n+2}}}-\frac{B(\frac{\ell_3}{2^{n}})-B(\frac{\ell_3+1}{2^{n}})}{2} - \frac{X_{(2\ell_3+1)2^{-(n+1)}}}{\sqrt{2^{n+2}}}
$$

\Bin we get that the expectation of the product of each term by $Y2$ is zero. For example, the expectations of the product of $B((\ell_2)2^{-n})$ by $Y2$ is

\begin{eqnarray*}
\frac{1}{4}\biggr(\frac{\ell_2}{2^{n}}+\frac{\ell_2}{2^{n}}+0-\frac{\ell_2}{2^{n}}-\frac{\ell_2}{2^{n}}-0\biggr)=0.
\end{eqnarray*}

\Bin The computations are the same for the other six cases since the ordering of the $\ell_i$'s follows the strict ordering of $r_i$. For example $2\ell_1 < 2\ell_2+1$ and $2\ell_1<2\ell_2$ both imply $\ell_1\leq \ell_2$.\\

\Ni \textit{(c) Checking (IS)}. For $(r_1,r_2) \in (D_{n+1}^{\star})^2$, such that $0\leq r_1<r_2$. Since the process is centered and Gaussian, checking (IS) reduces to checking that

$$
\mathbb{V}ar(B(r_2)-B(r_1))=r_2-r_1.
$$

\Bin But, here again, for each $i \in \{1,2\}$, we may have either $r_i=(2\ell_i)2^{-(n+1)}$ or $r_i=(2\ell_i+1)2^{-(n+1)}$. So we have four cases to check. But we may exclude the case in which $r_i=(2\ell_i)2^{-(n+1)}$ for all $i \in \{1,2\}$ since this is covered by the induction hypothesis. Here, we only check the case
$r_i=(2\ell_i+1)2^{n+1}$ for all $i \in \{1,2\}$ for peace of mind. The treatment of the two other cases are very similar in we mentioned in .

\Bin We have that  $B(r_2)-B(r_1)$ is $Z_1-Z2+Z3$ with

\begin{eqnarray*} 
Z_1&=&\frac{B(\frac{\ell_2}{2^{n}})+B(\frac{\ell_2+1}{2^{n}})}{2} \\
Z_2&=&\frac{B(\frac{\ell_1}{2^{n}})+B(\frac{\ell_1+1}{2^{n}})}{2}\\
Z_3&=&\frac{X_{(2\ell_2+1)2^{-n}}}{\sqrt{2^{n+2}}}-\frac{X_{(2\ell_1+1)2^{-n}}}{\sqrt{2^{n+2}}}.
\end{eqnarray*}

\Bin The variance of $Z_3$ is $2\times 2^{-(n+2)}$ and its covariance with $Z_1$ and $Z_2$ are null. The covariance of $Z1$ and $Z2$ is

$$
\frac{1}{4}\left(\frac{\ell_1}{2^{n}}+\frac{\ell_1}{2^{n}}+\frac{\ell_1}{2^{n}}+\frac{\ell_1+1}{2^{n}}\right)=\frac{4\ell_1+1}{2^{n+2}}.
$$

\Bin For each $i \in \{1,2\}$, the variances of $Z_i$ is

$$
\frac{1}{4}\left(\frac{\ell_i}{2^{n}} + \frac{\ell_i+1}{2^{n}}+2\frac{\ell_i}{2^{n}}\right)=\frac{1}{4}\left(\frac{3\ell_i}{2^{n}} + \frac{\ell_i+1}{2^{n}}\right).
$$

\Bin Finally, the variance of $B(r_2)-B(r_1)$ is

$$
\biggr(\frac{3\ell_2}{2^{n+2}} + \frac{\ell_2+1}{2^{n+2}}\biggr)+\biggr(\frac{3\ell_1}{2^{n+2}} + \frac{\ell_1+1}{2^{n+2}}\biggr)
-2\biggr(\frac{4\ell_1+1}{2^{n+2}}\biggr) + \biggr(\frac{2}{2^{n+2}}\biggr),
$$

\Bin which is

$$
\frac{4\ell_2-4\ell_1+2}{2^{n+2}}=\frac{2(\ell_2-\ell_1)+1}{2^{n+1}}=r_2-r_1.
$$

\Bin The proof is finished.\\

\Bin From there, the arguments in \cite{billinsgleypm} below allow the extension of the Broanian motion of $D$ on the space time $\mathbb{R}_+$ in a continuous version.

\section{Standard Brownian Motion} \label{cmb_03_sec_01}

\noindent Without any knowledge of the almost-sure continuity of the paths of the Brownian motion constructed by the \textit{KET}, we want to build a second version
with \textit{a.s} paths by using a weak convergence argument. Let us proceed into three steps.\\

\noindent \textbf{Step 1 : Construction of a co-null-set where $B$ has locally uniformly continuous paths on bounded intervals $[0,t]$ ($t>0$) of Dyadic numbers}.\\
 
\noindent  First let us consider the class all dyadic numbers in $\mathbb{R}_+$,

$$
D=\bigcup_{k\geq 0, \ n\geq 0} \left\{\frac{k}{2^n}\right\},
$$

\bigskip \noindent which is dense in $\mathbb{R}_+$. For all $(k,n)\in \mathbb{N}^2$, we define

$$
I_{k,n}=\left[ \frac{k}{2^n}, \frac{k+2}{2^n} \right],
$$

$$
M_{k,n}=\sup_{r \in I_{k,n} \cap D} \left|B\left(r\right)-B\left(\frac{k}{2^n}\right)\right| 
$$

\noindent and

$$
M_n=\sup_{k \leq (n+1)2^n} M_{k,n}. 
$$

\noindent We want to prove the : \\

\noindent \textbf{Fact 1}. $\mathbb{P}(M_n>n^{-1}, \ i.o.)=0$. $\Diamond$ \label{fact1P}\\

\noindent \textbf{Proof of Fact 1}. Let us fix $k \in \mathbb{N}$ and $n\geq 0$. Since $D \cap I_{k,n}$ is countable, we may write it as $\{r_1,r_2,....\}$ and we have

\begin{equation}
M_{k,n,p}=\sup_{1\leq i \leq p} |B\left(r_i\right)-B\left(\frac{k}{2^n}\right)| \nearrow  M_{k,n} \ as \ p \nearrow +\infty. \label{ineqPart}
\end{equation}

\Bin Now, for each $p$ fixed, we may set $r_0=k2^{-n}$ and form

$$
X_i=\biggr(B\left(r_i\right)-B\left(\frac{k}{2^n}\right)\biggr)- \biggr(B\left(r_{i-1}\right)-B\left(\frac{k}{2^n}\right)\biggr).
$$

\bigskip \noindent  The random variables $X_i$ are independent, by the independence of the increments of the Brownian motion. We also have

$$
B\left(r_i\right)-B\left(\frac{k}{2^n}\right)=X_1+...+X_i.
$$

\bigskip \noindent Hence, we may use Etemadi's Inequality (see \cite{ips-mfpt-ang}, Chapter 6, Inequality 18), that is, for any $\alpha>0$,

$$
\mathbb{P}\left(\max_{1\leq k \leq n} |X_1+\dots+X_k| \geq 3\alpha\right) \leq 3 \mathbb{P}\left(|X_1+\dots+X_n|\geq \alpha \right).
$$

\bigskip \noindent When applied here, we have for any $\alpha>0$,

$$
\mathbb{P}(M_{k,n,p} \geq 3 \alpha)\leq 3 \mathbb{P}\left(\left|B\left(r_p\right)-B\left(\frac{k}{2^n}\right) \right| \geq \alpha\right)
$$

\bigskip \noindent and by using Markov's inequality we get

\begin{eqnarray*}
\mathbb{P}(M_{k,n,p} \geq 3 \alpha) &\leq& 3 \mathbb{P}\left(\left|B\left(r_p\right)-B\left(\frac{k}{2^n}\right)\right| \geq \alpha\right)\\
&=& 3 \mathbb{P}\left(|N(0,1)| \geq \alpha (r_p-\frac{k}{2^n})^{-1/2}\right)\\
&=& 3 \mathbb{P}\left(|N(0,1)|^4 \geq \alpha^4 (r_p-\frac{k}{2^n})^{-2}\right)\\
&\leq& 9 \alpha^{-4} \left(r_p-\frac{k}{2^n}\right)^{2}\\
&\leq& 36 \alpha^{-4} 2^{-2n}.\\
\end{eqnarray*}

\noindent By letting $p\nearrow +\infty$, by using the continuity of the probability in Formula \ref{ineqPart}, we get

\begin{eqnarray*}
\mathbb{P}(M_{k,n} \geq 3 \alpha) &\leq& 36 \alpha^{-4} 2^{-2n}.
\end{eqnarray*}

\bigskip \noindent Finally, we have

\begin{equation*}
\mathbb{P}(M_{n} \geq 3\alpha) \leq \sum_{k=1}^{(n+1)2^n} \mathbb{P}(M_{k,n} \geq 3 \alpha) \leq 36 (n+1)\alpha^{-4} 2^{-n}.
\end{equation*}

\bigskip \noindent For $\alpha=n^{-1}/3$, we get

\begin{eqnarray*}
\mathbb{P}(M_{n} \geq n^{-1})  \leq 2916 (n+1)n^{4} 2^{-n}.
\end{eqnarray*}

\bigskip \noindent The convergence of the series of general term $\mathbb{P}(M_{n} \geq n^{-1})$ finishes the proof. $\square$\\

\noindent \textbf{Step 2 : Construction of continuous version of $B$}.\\

\noindent Now, let us denote $N=\{M_n > n^{-1}, \ i.o.\}$. Fix $\omega \in N^c$, $t \in \mathbb{R}_+, \ t>0$, and $\varepsilon>0$. So $M_n(\omega) \leq n^{-1}$ when $n$ is larger than some $n_0(\omega)$. And we may and do pick one those $n$'s such that

$$
t<n \ \ and \ \ 2n^{-1}\leq \varepsilon.
$$

\bigskip \noindent Let us fix $\delta=2^{-n}$. We have

$$
[0,t]=\{0\} + \biggr(\sum_{0\leq k \leq K(t)-1} \left]\frac{k}{n}, \frac{k+1}{n}\right]\biggr)+ \left]\frac{K(t)}{n}, t\right], \ for \ \frac{K(t)}{n} \leq t <\frac{K(t)+1}{n}. 
$$

\bigskip \noindent It is clear that if $r_1$ and $r_2$ are two elements of $[0,t]$ such that $|r_1-r_2|<2^{-n}$, and since the intervals which partition $[0,t]$ are at most of length $2^{-n}$, we deduce that $r_1$ and $r_2$ are in one interval of the decomposition or, in the worse case, in two adjacent intervals, meaning that they are in some $I_{n,k}$, where 

$$
\frac{k}{n} \leq \frac{K(t)+1}{n}\leq t+n^{-1},
$$

\bigskip \noindent which implies

$$
k \leq (n+1)2^{n}.
$$

\bigskip \noindent It follows that

\begin{eqnarray*}
|B(r_1,\omega)-B(r_2,\omega)|&=&\left|\biggr(B\left(r_1\right)-B\left(\frac{k}{2^n}\right)\biggr)- \biggr(B\left(r_{r_2}\right)-B\left(\frac{k}{2^n}\right)\biggr)\right|\\
&\leq& 2M_n(\omega)\leq 2n^{-1}\leq \varepsilon.
\end{eqnarray*}

\bigskip \noindent We just proved that $D \ni r \rightarrow B(\omega,r)$ is uniformly continuous on each $[0,t]$, $t>0$, outside the null set $N$.\\

\noindent \textbf{Step 3 : Use of Weak Convergence Theory to conclude}.\\

\noindent For each $t\geq 0$, for $\omega \in N^c$, the sequence $\{B(\omega,r), \ r \in D \cap ]t, t+1]\}$, as $r \searrow t$, has the Cauchy property and hence, we may define

\begin{equation}
\tilde{B}(t,\omega)=\lim_{D \ni r \searrow t} B(r,\omega). \label{limitBM}
\end{equation}

\bigskip \noindent By continuity of : $D \ni r \rightarrow B(\omega,r)$, we have $\tilde{B}=B$ on $N^c$. That $\tilde{B}(\circ,\omega)$ is continuous for $\omega \in N^c$ is direct. Indeed let   $(t_1,t_2) \in \mathbb{R}_+^2$ such that $0\leq t_1\leq t_2\leq t$. Let $\varepsilon>0$. By the uniform continuity of $B(r,\omega)$ on 
$[0,t]\cap D$, there exists $\delta>0$ such that $|B(r,\omega)-B(s,\omega)|<\varepsilon/3$ whenever $(r,s)\in ([0,t]\cap D)^2$ and $|r-s|<\delta$. Now suppose that $0\leq t_2-t_1\leq \delta/2$. Also, from the density of the dyadic numbers in $\mathbb{R}_+$ and by the limit in Formula \eqref{limitBM}, $0\leq t_2-t_1<\delta$, we may find $r_1$ and $r_2$ in $D$ such that

$$
t_2< r_2, \ t_1<r_1,  \ \ 0\leq r_2-r_1\leq 2\delta \ \ and \ \ |B(r_i,\omega)-\tilde{B}(t_i,\omega)|\leq \varepsilon/3, \ 1\leq i \leq 2.
$$

\Bin Hence $|\tilde{B}(t_1,\omega)-\tilde{B}(t_2,\omega)|\leq \varepsilon$. So $\tilde{B}(\circ,\omega)$ is continous. Furthermore, $\tilde{B}$ and $B$ have the finite distribution probability laws. Indeed, let $0=t_0<t_1<...<t_k$, $k\geq 1$. Since

$$
(B(r_1), B(r_2),...,B(r_k)) \rightarrow (\tilde{B}(t_1),\tilde{B}(t_2),...,\tilde{B}(t_k)) \ as \ r_1\searrow t_1, r_2\searrow t_2, ...,r_k\searrow t_k,
$$

\bigskip \noindent outside the null-set $N$, we get, by the comparison Theorem 3 in \cite{ips-wcia-ang}, Chapter 2, page 59 [of convergence in probability and convergence in law of random vectors defined on the same probability space] then

$$
(B(r_1), B(r_2),...,B(r_k)) \rightsquigarrow (\tilde{B}(t_1),\tilde{B}(t_2),...,\tilde{B}(t_k)) \ as \ r_1\searrow t_1, r_1\searrow t_1, ...,r_k\searrow t_k.
$$

\bigskip \noindent By the Portmanteau Theorem 4 in \cite{ips-wcia-ang} (Chapter 2, page 61), the cumulative distribution functions \textit{cdf}'s converge, that is for any continuity point $x=(x_1,...,x_k)$ of the \textit{cdf} $F_{\tilde{B}(t_1),\tilde{B}(t_2),...,\tilde{B}(t_k)}$, we have
 
\begin{equation}
F_{(\tilde{B}(t_1),\tilde{B}(t_2),...,\tilde{B}(t_k))}(x_1, x_2,...,x_k)=\lim_{r_1\searrow t_1,...,r_k\searrow t_k} F_{(B(r_1),B(r_2),...,B(r_k))}(x_1, x_2,...,x_k). \label{limitBM1}
\end{equation}

\bigskip \noindent Now, for any $(s_1,...,s_k) \in \mathbb{R}^k$, 

\begin{eqnarray*}
f_{(B(r_1),...,B(r_k))}(s_1,...,s_k)&=&\prod_{j=1}^{k} \frac{1}{\sqrt{2\pi(r_j-r_{j-1})}} \exp\left(-\frac{(s_j-s_{j-1})^2}{2(r_j-r_{j-1})}\right)\\
&\rightarrow&\prod_{j=1}^{k} \frac{1}{\sqrt{2\pi(t_j-t_{t-1})}} \exp\left(-\frac{(s_j-s_{j-1})^2}{2(t_j-t_{j-1})}\right)\\
&=&f_{(B(t_1),...,B(t_k))}(t_1,\cdots,t_k).
\end{eqnarray*}

\Bin By Scheff\'e's Theorem (proposition 13  in \cite{ips-wcia-ang}, Chapter 2, page 67), we have, for any continuity point $(x_1, x_2,...,x_k)$ of $F_{B(t_1), \ B(t_2),..., \ B(t_k)}$,

\begin{equation}
F_{(B(r_1),B(r_2),...,B(r_k))}(x_1, x_2,...,x_k)\rightarrow F_{(B(t_1), \ B(t_2),..., \ B(t_k))}(x_1, x_2,...,x_k).\label{limitBM2}
\end{equation}

\bigskip \noindent By combining Formulas \eqref{limitBM1} and \eqref{limitBM2}, we have, for any continuity point $(x_1, x_2,...,x_k)$ of $F_{(B(t_1), \ B(t_2),\cdots, \ B(t_k))}$ and for any continuity point $(x_1, x_2,...,x_k)$ of $F_{(\tilde{B}(t_1), \ \tilde{B}(t_2),..., \ \tilde{B}(t_k))}$, 

$$
F_{(\tilde{B}(t_1),\tilde{B}(t_2),...,\tilde{B}(t_k))}(x_1, x_2,...,x_k)=F_{(B(t_1), \ B(t_2),..., \ B(t_k))}(x_1, x_2,...,x_k).
$$

\Bin We conclude that

$$
F_{(B(t_1),B(t_2),...,B(t_k))}=F_{(\tilde{B}(t_1),\tilde{B}(t_2),...,\tilde{B}(t_k))}.
$$

\bigskip \noindent Hence, The \textit{a.s.} continuous paths stochastic process $\tilde{B}$ and $B$ have the same finite distribution laws.\\

\noindent \textbf{Conclusion}. The stochastic process $\tilde{B}$ has \textit{a.s.} continuous paths and has the Brownian motion probability laws.\\

\section{Conclusion} \label{cmb_02_sec_01}

\Ni The full an detailed construction of the Standard Brownian motion is complete. The only external tool used in the paper is the proof of the Etamadi's inequality. to make the paper more self-readable, we add that proof in the appendix.\\

\newpage

\Bin \textbf{Appendix : Etemadi's Inequality}. Let  $X_{1}, \cdots, \leq X_{n}$ be $n$ independent real-valued random variables such that the partial sums $S_{k}=X_{1}+...+X_{k}$, $1\leq k\leq n$, are definied. Then for any $\alpha \geq 0$, we have

\begin{equation*}
\mathbb{P}\left(\max_{1\leq k\leq n}\left\vert S_{k}\right\vert \geq 3\alpha \right)\leq 3 \max_{1\leq k\leq n}\mathbb{P}\left(\left\vert S_{k}\right\vert \geq \alpha\right). \ \Diamond
\end{equation*}

\Bin \textbf{Proof}. The formula is obvious for $n=1$. Let $n\geq 2$. As usual, denote $B_1=(|X_1|\leq 3\alpha)$, $B_k=(|S_1|<3\alpha, \cdots , |S_{k-1}|<3\alpha, |S_k|\geq  3\alpha)$, $k\geq 2$. By decomposing $(\max_{1\leq j \leq n} |S_j| \geq 3\alpha)$ over the partition 
$$
(|S_n| \geq \alpha)+(|S_n| < \alpha)=\Omega,
$$

\Bin we have

$$
(\max_{1\leq j \leq n} |S_j| \geq 3\alpha) \subset (|S_n| \geq \alpha) \cup (|S_n|<\alpha, \max_{1\leq j \leq n} |S_j| \geq 3\alpha)
$$

\Bin And by the principle of the construction of the $B_j$,

$$
(\max_{1\leq j \leq n} |S_j| \geq 3\alpha)=\sum_{1\leq j \ n} B_j
$$ 

\Bin and hence

$$
(\max_{1\leq j \leq n} |S_j| \geq 3\alpha) \subset (|S_n| \geq \alpha) \cup \sum_{1\leq j \ n-1} (|S_n| < 3\alpha) \cap B_j
$$

\Bin where the summation is restricted to $j \in \{1,...,n-1\}$ since the event $(|S_n| < 3\alpha) \cap B_n$ is empty. Further, on $(|S_n| < \alpha) \cup B_j$, we have $(|S_n| < \alpha)$ and $(|S_j| < 3\alpha)$ and the second triangle inequality $|S_n-S_j|\geq |S_j|-|S_n|\geq 3\alpha-\alpha=2\alpha$, that is

$$
(|S_n| < \alpha) \cap B_j \subset B_j \cap (|S_n|<\alpha) \cap (|S_n-S_j|\geq 2\alpha) \subset B_j \cap (|S_n-S_j|\geq 2\alpha).
$$

\Bin Now, we remind that and $B_j$ and $S_n-S_j$ are independent. Translating all this into probabilities gives

\begin{eqnarray*}
\mathbb{P}(\max_{1\leq j \leq n} |S_j| \geq 3\alpha)&\leq&  \mathbb{P}(|S_n| \geq \alpha) + \sum_{1\leq j \ n-1} \mathbb{P}(B_j) \mathbb{P}(|S_n-S_j|\geq 2\alpha)\\
&\leq&  \mathbb{P}(|S_n| \geq \alpha) + \sum_{1\leq j \ n-1} \mathbb{P}(B_j) \biggr(\mathbb{P}(|S_n|\geq \alpha)+\mathbb{P}(|S_j|\geq \alpha)\biggr).\\
\end{eqnarray*} 

\Bin But $(|S_n|\geq 2\alpha)$, $(|S_n| \geq \alpha)$ and $(|S_j|\geq 2\alpha)$ are subsets of 
$$
(\max_{1\leq j \leq n} |S_j| \geq \alpha)
$$ 

\Bin and hence, we may conclude that

\begin{eqnarray*}
\mathbb{P}\left(\max_{1\leq j \leq n} |S_j| \geq 3\alpha\right)&\leq&  \mathbb{P}\left(\max_{1\leq j \leq n} |S_j| \geq \alpha\right) \left(1+ 2 \sum_{1\leq j \leq n}\mathbb{P}(B_j)\right)\\
&\leq&  \mathbb{P}\left(\max_{1\leq j \leq n} |S_j| \geq \alpha\right) \left(1+ 2 \mathbb{P}\left(\sum_{1\leq j n}B_j\right)\right)\\
&\leq&  3\mathbb{P}(\max_{1\leq j \leq n} |S_j| \geq \alpha). \ \square
\end{eqnarray*} 
\end{document}